\documentclass[aip,cha,preprint,numerical, nofootinbib, groupedaddress]{revtex4-1}
\usepackage{amsmath, amsthm, amssymb}
\usepackage{graphicx}
\usepackage{listings}
\usepackage{subfigure}

\newcommand{\Eq}[1]{(\ref{eq:#1})}
\newcommand{\Fig}[1]{Fig.~\ref{fig:#1}}

\newcommand{\NN}{{\cal N}}

\newtheorem*{defn}{Definition}

\begin{document}
	
	\title{Iterated Function System Models in Data Analysis: Detection and Separation}	
	\author{Zachary Alexander}
	\email[Electronic Address: ]{alexanz@colorado.edu}	
	\author{Elizabeth Bradley}
	\email[Electronic Address: ]{lizb@cs.colorado.edu}	
	\author{Joshua Garland}
	\email[Electronic Address: ]{joshua.garland@colorado.edu}
	\author{James D. Meiss}
	\email[Electronic Address: ]{jdm@colorado.edu}	
	\affiliation{University of Colorado}

\begin{abstract}
We investigate the use of iterated function system (IFS) models for data analysis. An IFS is a discrete dynamical system in which each time step corresponds to the application of one of a finite collection of maps. The maps, which represent distinct dynamical regimes, may act in some pre-determined sequence or may be applied in random order.
An algorithm is developed to detect the sequence of regime switches under the assumption of continuity. This method is tested on a simple IFS and applied to an experimental computer performance data set. This methodology has a wide range of potential uses: from change-point detection in time-series data to the field of digital communications. 
\end{abstract}

\maketitle

\begin{quotation} 
	An IFS is a discrete dynamical system in which one member of a finite
	collection of maps acts at each time step.  The maps, which represent
	distinct dynamical regimes of the overall system, may act randomly or
	in some pre-determined sequence.  This is a useful framework for
	understanding the dynamics of a wide class of interesting systems,
	ranging from digital communication channels to the human brain.  In
	this paper, we first review the IFS framework and then present an
	algorithm that leverages the associated properties to segment the
	signal.  Working under the assumption that each of the maps is
	continuous, this algorithm uses topology to detect the different
	components in the output data.  The main idea behind this approach is
	that nearby state-space points can evolve in different ways, depending
	on the dynamical regime of the IFS; the main challenge is that the
	components may overlap, causing their trajectories to locally
	coincide.  We demonstrate this algorithm on two examples, a H\'enon IFS
	and an experimental computer performance data set.  In both cases, we
	were able to segment the signals effectively, identifying not only the
	times at which the dynamics switches between regimes, but also the
	number and forms of the deterministic components themselves.

\end{quotation}

\section{Introduction} \label{sec:introduction}

Any approach to time series analysis begins with the question: is the
data stochastic or deterministic \cite{theiler1992testing, sugihara1990nonlinear}? Often, the answer may be ``both'': the data could be generated by a deterministic system with a noisy component, perhaps due to measurement or computer round-off error. In this article, we propose an alternative possibility: the data could be generated by a sequence of deterministic dynamical systems selected by a switching process that itself could be deterministic or stochastic, i.e., by an \emph{Iterated Function System} (IFS). For a detailed review of IFS dynamics, see Diaconis and Freedman\cite{diaconis1999iterated}. If this were the case, then a useful goal is to identify the times at which switching between regimes occurs, as well as the number and forms of the deterministic components themselves.  Under the assumption that each deterministic system is continuous, we use topology to detect and separate the components of the IFS that are present in the output data. The main idea behind this approach is that the nearby state-space points can evolve in different ways, depending on the dynamical state of the IFS. Such a model has some relation to the determination of states in a hidden Markov model; however, hidden Markov models are typically discrete and stochastic---not continuous and deterministic.  A primary
challenge in this problem is that overlap between the images of distinct regime functions could cause their trajectories to locally coincide. The use of IFS models for physical systems is not new; for example, Broomhead, et al.\cite{broomhead2004iterated} used an IFS to
model digital communication channels. In the current paper, we provide new tools to extract IFS models from experimental data and to determine the sequence of switching between regimes in the IFS.

We believe the method proposed here will prove useful in a number of applications. For example, detection and separation of IFS components is closely related to the statistical problem of \emph{event} or \emph{change-point detection} \cite{kawahara2009change} where time-series data is assumed to come from a statistical distribution that changes suddenly. Applications where change-point detection plays a role include fraud detection in cellular systems, intrusion detection in computer networks, irregular-motion detection in computer vision, and fault detection in engineering systems, among many others \cite{kawahara2009change}. Our underlying hypothesis is different than in the field of statistical change-point detection---we assume that each regime is deterministic. For example, though change-point detection has been successfully applied to determine brain states from EEG data \cite{kaplan2000application}, EEGs have also been shown to exhibit properties of low-dimensional chaos \cite{klonowski1999quantitative}. Indeed, low-dimensional dynamics
occurs in diverse areas including physiology, ecology, and economics \cite{glass1993time,turchin1992complex, chen1988empirical}. We expect that the separation technique outlined below could be used to produce more-accurate models of regime shifts and the effects of rapid parameter changes that occur, e.g. in the onset of seizures, natural disasters, or the bursting of economic bubbles.


\section{Detection and Separation} \label{sec:iterated:function:systems}

Given a time series that corresponds to measurements of a dynamical system, our goal is to develop a technique that will detect whether the series is generated by an Iterated Function System (IFS) and to distinguish its components. Formally, an IFS is a discrete-time dynamical system that consists of a finite set of maps $\{f_0,\ldots, f_n,\ldots, f_{N-1}\}$ on a state space $X$. A trajectory of the IFS is a sequence of state-space points, $\{x_0,\ldots, x_t,x_{t+1},\ldots\}$, together with a \emph{regime} sequence $\{n_0,\ldots, n_t,n_{t+1},\ldots\}$ with $n_t \in \{0,1,\ldots,{N-1}\}$, such that
\[
	x_{t+1} = f_{n_t}(x_t) \;, \quad \forall t \in \mathbb{N} .
\]
Without loss of generality, we may assume that each map occurs at least once in the regime sequence, since otherwise the missing maps could be eliminated. 

In the standard study of IFS dynamics, the regime sequence is often taken to be a realization of some random process \cite{barsnley1988fractals, falconer2003fractal}; however, we only assume that we have access to a single trajectory that is generated by a particular realization. Consequently, the selection rule for the regime sequence is immaterial; indeed, it could just as well be a discrete, deterministic dynamical system. The standard theory, in addition, often requires that each $f_n$ is a contraction mapping, 
in which case the IFS is \emph{hyperbolic} and has a unique attractor $A$ that is invariant in the sense that  $A = \bigcup_{n=0}^{N-1} f_n(A)$. We do not not need this assumption, and only assume that the trajectory lies in some bounded region of $X$. 

We will assume that the time series corresponds to $T$ measurements on a particular state-space sequence, 
\begin{equation}\label{eq:Gamma}
	\Gamma = \{x_0,x_1,\ldots, x_{T-1}\};
\end{equation}
but that the regime sequence is unmeasurable or hidden. For example, one may be able to 
measure the position of a forced pendulum at a sequence of times, but the pendulum may have 
a sealed brake mechanism that sets a friction coefficient and that is controlled externally to the experiment. Measurement of $\Gamma$ also implicitly includes that of its associated shift map
\begin{equation}\label{eq:sigma}
	\sigma(x_t) = x_{t+1}.
\end{equation} 

It is often the case that a time series corresponds to a limited measurement, perhaps of  one variable from a multi-dimensional dynamical system. In this case, the first step is to use delay-coordinate embedding to construct, as much as is possible, a topologically faithful image of the orbit a reconstructed state-space \cite{Sauer:1991lr,Takens:1981uq}. We suppose that \Eq{Gamma} is this embedded time series. 

The fundamental goal in this paper is \emph{detection and separation}: to detect 
if $\Gamma$ is a trajectory of an IFS and to separate the regimes by 
recovering the sequence $\{n_t\}$.
This problem is relatively straightforward when $\Gamma$
is a subset of some \emph{non-overlapping} region of the IFS, i.e., a region $R$ such that $f_i(R) \bigcap f_j(R) = \emptyset$ for all $i \neq j$.
In this paper, we address a more-general situation in which $\Gamma$ could be sampled from an overlapping region of the IFS.

A fundamental requirement for our separation method is that the maps $f_n$ are continuous---a reasonable assumption for the vast majority of physical systems. In particular, the image of a connected set under each $f_n$ must be connected. Since for finite data sets, the notion of connectivity makes no sense, we will instead use $\epsilon$-connectivity under the assumption that $X$ is a metric space with distance $d(x,y)$. 

\begin{defn}[$\epsilon$-connected \cite{robins2000computing}] A set $\Omega \subset X$ is \emph{$\epsilon$-connected} if there exists an $\epsilon$-chain connecting the points in $\Omega$, i.e., for each pair of points $x,y \in \Omega$ there exists a sequence $\{z_0,\ldots,z_k\}\subset \Omega$ such that $x = z_0$, $y=z_k$, and $d(z_j,z_{j+1}) < \epsilon$ for $0 \le j \le k -1$.
\end{defn}

\noindent
Let $\NN_k(x_t)$ denote the set of $k$ points consisting of $x_t$ and its $(k-1)$-nearest neighbors in $\Gamma$. For each such set there will be a $\delta$ such that $\NN_k(x_t)$ is $\delta$-connected.

The idea of our algorithm is as follows. For each $\epsilon$ that is not too small, there must be a $k>1$ such that the image of $\NN_k(x_t)$ under a single map will be $\epsilon$-connected. Indeed, continuity implies there is a $\delta$ such that a $\delta$-connected set has an $\epsilon$-connected image.  For a given $\epsilon$, the minimal $\delta$ will be determined by the maximal distortion of the map. For the algorithm to work, the set $\Gamma$ must be dense enough so that for this $\delta$, there are nearest neighbors, i.e., $k > 1$.

\begin{figure}[ht]
	\centering
	\includegraphics[scale=.35]{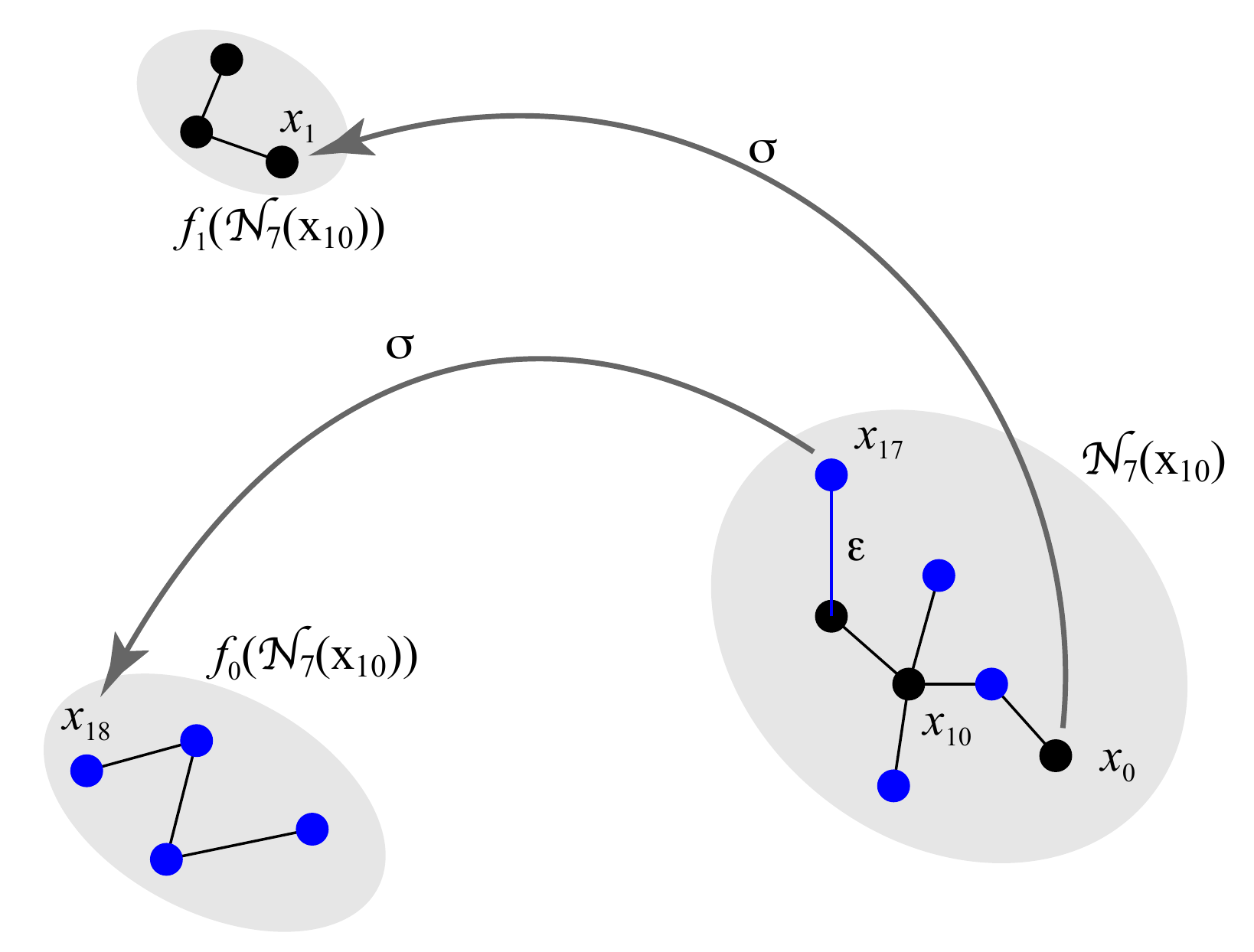}
	\caption{Sketch of the action of the shift map $\sigma$ on a $7$-nearest neighborhood of a point $x_{10} \subset \Gamma$ results in two $\epsilon$-connected components that can be identified as $f_0(\NN_7(x_{10}))$ and $f_1(\NN_7(x_{10}))$.
\label{fig:sepsketch}}
\end{figure}

If $\epsilon$ is chosen to reflect this maximal, single-map distortion, then whenever the time-shifted image, $\sigma(\NN_k(x_t))$, consists of a number of $\epsilon$-connected components, each component should reflect the action of a different $f_n$. This idea is expressed visually in \Fig{sepsketch}. Note that $\sigma(\NN_k(x_t))$ is NOT the same as $\NN_k(x_{t+1})$, the set of nearest neighbors to the image of $x_{t}$. 

To obtain reasonable results the parameter $\epsilon$ must be selected carefully as it will  determine the maximal number of nearest neighbors, $k$. The number, $N$, of regimes of the IFS is not more than the maximal number of components of $\sigma(\NN_k(x_t))$. However, since sparsely covered portions of the data set could result in spurious components, we will select $N$ to be the number of components in the bulk of the images $\sigma(\NN_k(x_t))$.

Given a time series $\Gamma$ that we suspect to be generated by an IFS, the detection and separation algorithm requires an appropriate value for $\epsilon$. Here we outline a possible algorithm.
\begin{itemize}
	\item (Detection) Determine a value for $\epsilon$ by computing histograms of the separations between $\NN_2(x_t)$ and $\sigma(\NN_2(x_t))$. If $\Gamma$ is sampled from a connected invariant set, then each of the nearest neighbor sets should be $\delta$-connected. If there is more than one regime, their images should be disconnected for some choice of $\epsilon$. The number of regimes $N$ is estimated to be the number of components in the majority of the $\sigma(\NN_k(x_t))$; this should be persistent over a range of $\epsilon$ and $k$ values.
   \item (Separation) Select a set of $K$-nearest-neighborhoods, $\{\Omega_j = \NN_K(x_{t_j})| j=0,1,\ldots J-1\}$, that overlap and cover $\Gamma$. Points are identified to be from a common regime if they lie in overlapping neighborhoods and their images lie in $\epsilon$-connected components.
\end{itemize}

In the next section, we illustrate this method on a simple example.

\section{Example: a H\'enon IFS} \label{sec:henon-like:ifs}
\begin{figure}[ht]
	\centering
	\includegraphics[scale=.40]{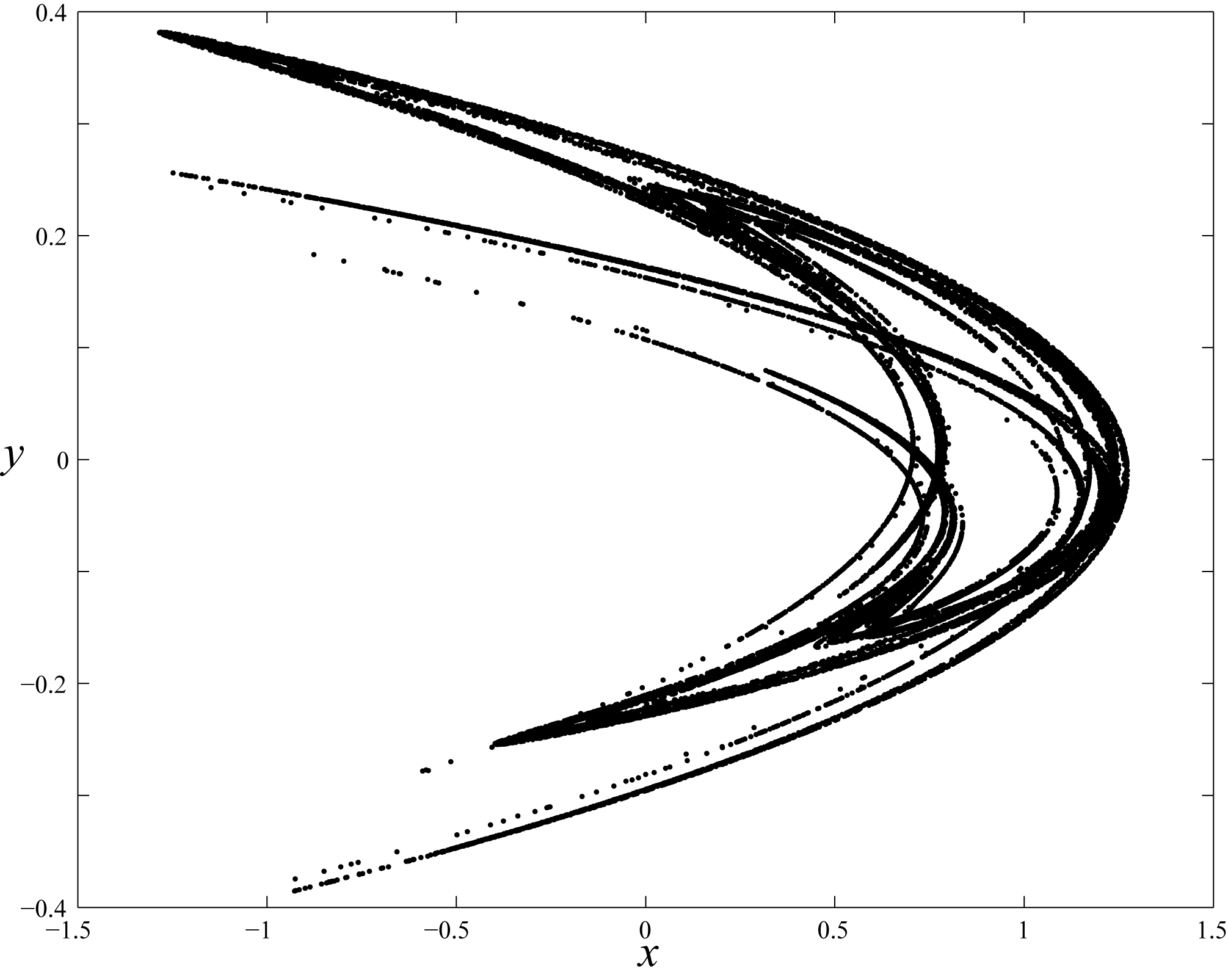}
	\caption{A trajectory of the IFS generated by \Eq{Henon} 
	with $T = 30,000$ points.  Here $n_t \in \{0,1\}$ are chosen with equal probability.  
	\label{fig:ifs:trajectory}}
\end{figure}

As a simple example, consider the IFS generated by the two quadratic, planar diffeomorphisms
\begin{equation}\begin{split}\label{eq:Henon}
	f_0 \left( x,y \right) &= 
	\left( y + 1 - 1.4x^2,\, 0.3x  \right), \\
	f_1 \left( x,y \right) &= 
	\left( \ y + 1 - 1.2(x-0.2)^2 ,\, -0.2x \right).
\end{split}\end{equation}
The map $f_0$ is H\'enon's quadratic map with the canonical choice of parameter values \cite{henon1976two}; the map $f_1$ is conjugate, via an affine change of coordinates, to H\'enon's map with parameters $(a,b) = (0.912,-0.2)$.  We generate a single trajectory of this IFS by using a Bernoulli process with equal probability to generate a sequence $n_t \in \{0,1\}$. A trajectory with $T= 30,000$ points, shown in \Fig{ifs:trajectory}, has the appearance of two overlapping H\'enon-like attractors. Note however, that since most points on $\Gamma$ are not iterated more than a couple of consecutive steps with the same map, $\Gamma$ is not just the union of the attractors of $f_0$ and $f_1$. Indeed the attractor for $f_1$ is simply a fixed point at $(0.63986,-0.12797)$.

\begin{figure}[ht]
	\centering
	\subfigure{
		\includegraphics[scale=.45]{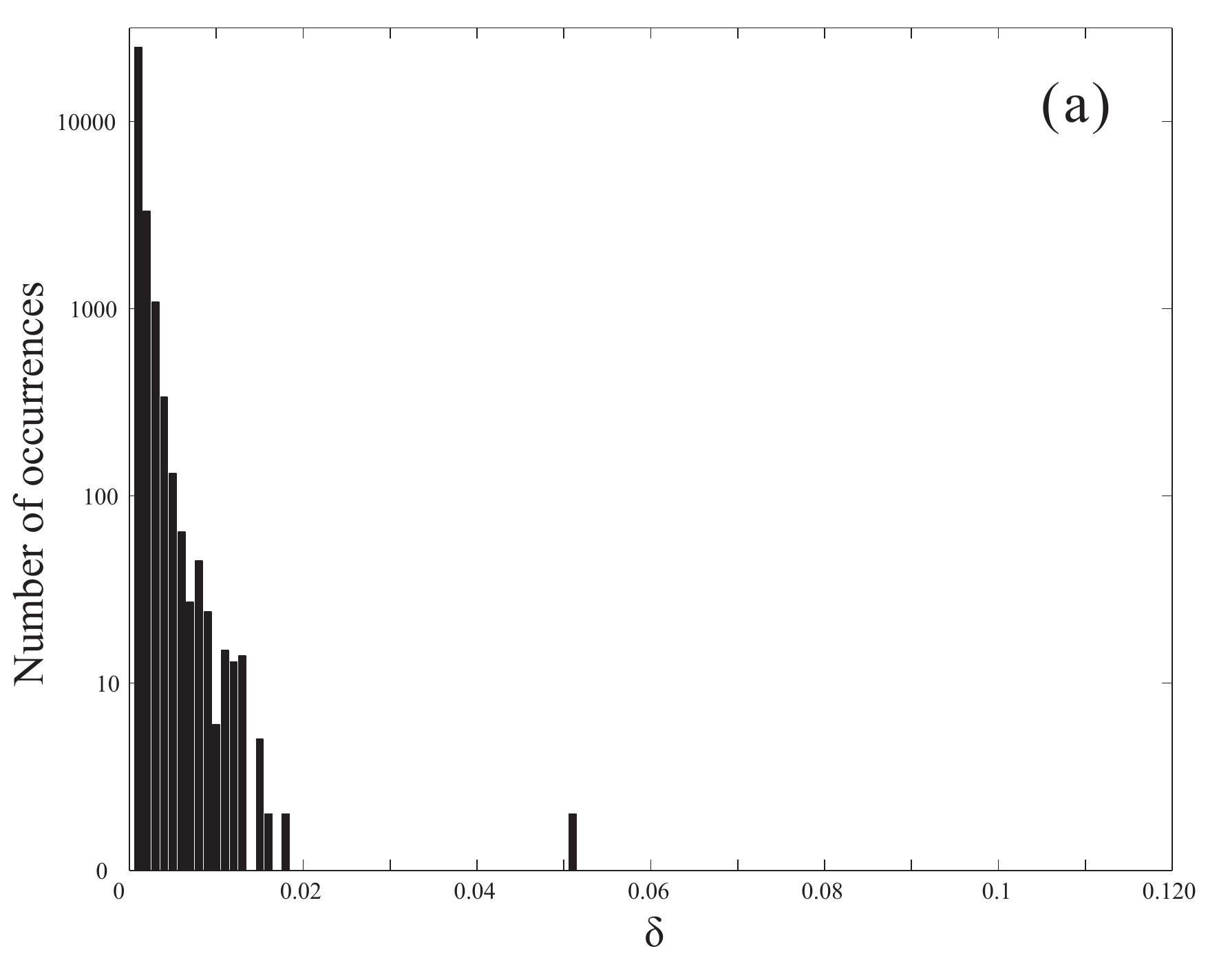}
	}
	\subfigure{
		\includegraphics[scale=.45]{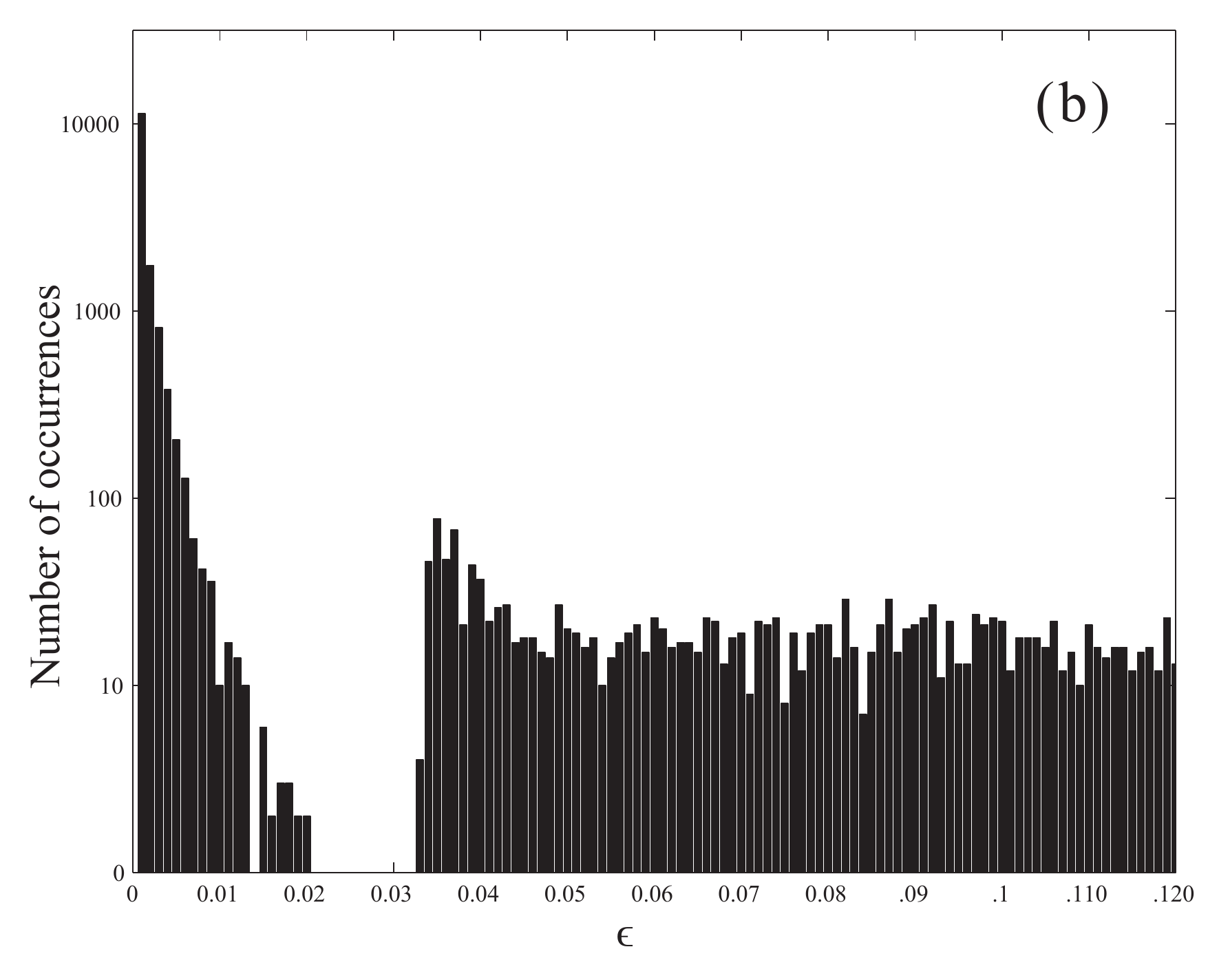}
	}
	\caption{Distance between each point on $\Gamma$ of \Fig{ifs:trajectory} and its nearest neighbor (a) and between the images of these two points (b). \label{fig:epsilonhist}}
\end{figure}

To recover the regime sequence from $\Gamma$ we must check for
$\epsilon$-disconnected images of $\delta$-connected components. As a first step to determine an appropriate value for $\epsilon$, we compute the distance between each point in $\Gamma$ and its nearest neighbor, i.e., the diameter of $\NN_2(x_t)$. This is shown in  panel (a) of \Fig{epsilonhist} as a histogram. Note that all but two points in $\Gamma$ have a nearest neighbor within $0.02$, and the vast majority within $0.002$.  Panel (b) of \Fig{epsilonhist} indicates how these distances grow upon iteration: it shows the distance between the iterates of each of these nearest neighbors, i.e., the diameter of $\sigma(\NN_2(x_t))$. There are now two distinct distributions separated by a gap $[0.02,0.032]$. This suggests that the dynamics underlying $\Gamma$ is discontinuous, and that a choice of $\epsilon$ in the gap may be appropriate. 

\begin{figure}[ht]
	\centering
	\includegraphics[scale=.5,trim=.75in 0 0 0]{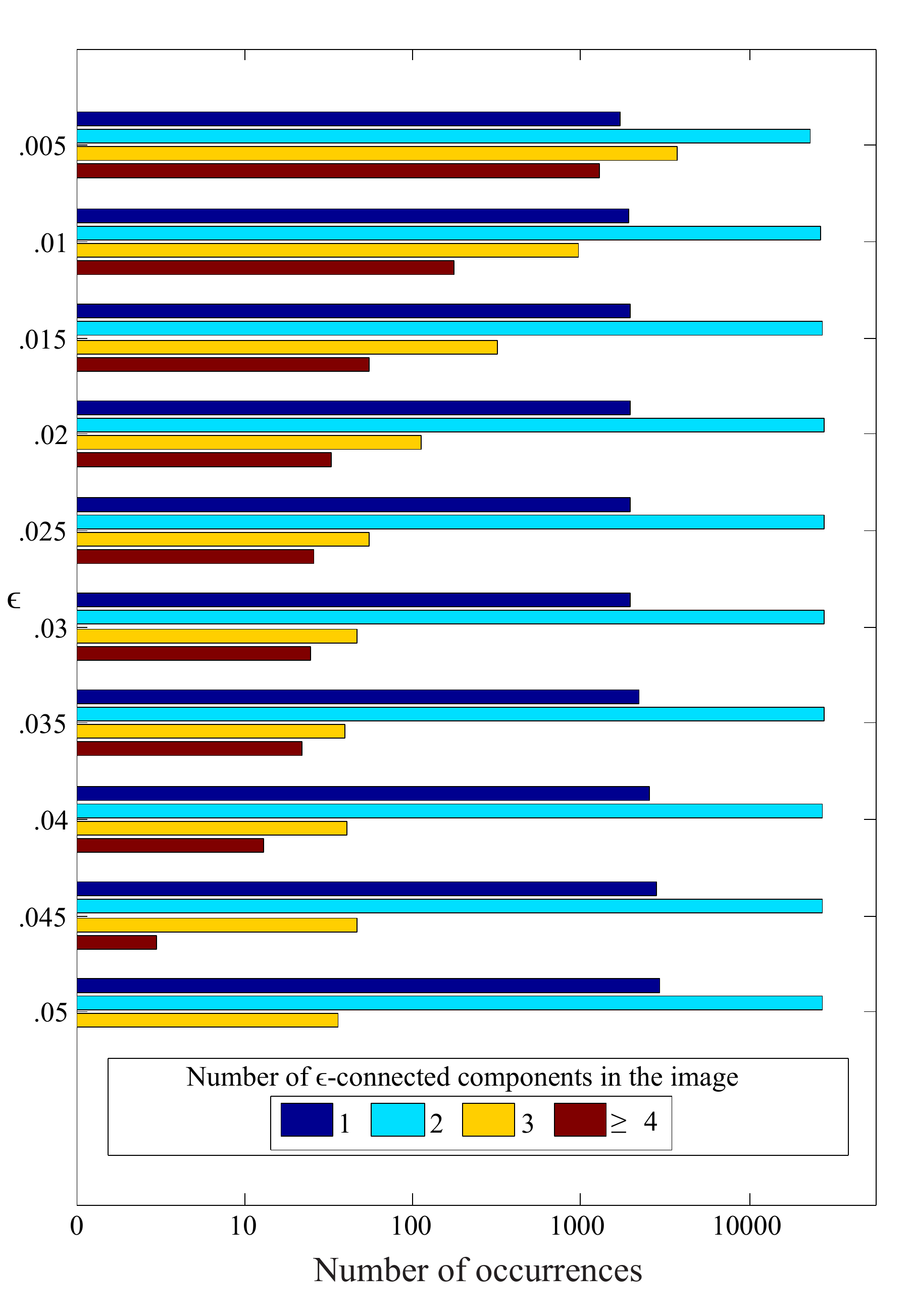}
	\caption{Detection of the number of components of the images of sets of $K=5$ nearest neighbors for the trajectory $\Gamma$ of \Fig{ifs:trajectory}. The histograms show the number of the images $\{\sigma(\NN_5(x_t))\;|\; 0\le t < 29,999\}$ that have $N$ $\epsilon$-connected components for various $\epsilon$.
\label{fig:henonhist}}
\end{figure}

Suppose that we did not know that the IFS \Eq{Henon} had two regimes---that only the trajectory $\Gamma$ of \Fig{ifs:trajectory} was available.
To detect the number of regimes, we look at the number of $\epsilon$-components in the image of the sets of five nearest-neighbors, $\NN_5(x_t)$. 
Histograms of the number of $\epsilon$-components of $\sigma(\NN_5(x_t))$ are shown in \Fig
{henonhist} as $\epsilon$ varies from $0.005$ to $0.05$. The vast majority of these neighborhoods split into at most two $\epsilon$-components. When $\epsilon$ is as small as $0.005$ about 3\% split into four or more components and when $\epsilon \ge 0.02$, only 0.3\% split into three or more components. Note that with the equal probability rule that we used for \Eq{Henon}, the probability that all five points in $\NN_5(x_t)$ will be iterated with the same map is $\tfrac{2}{32}$, which is confirmed in \Fig{henonhist}, since about 6\% of the images have one $\epsilon$-component.

Thus in the detection phase of the algorithm, we confirm that the underlying dynamics has two regimes, $N = 2$ and obtain a reasonable choice, $\epsilon = 0.03$.


For the separation phase of the algorithm, we wish to classify which points on $\Gamma$ are 
images of which map. To do this we choose larger, overlapping neigbhorhoods that 
cover $\Gamma$ so that we can connect the subsets for each regime. To distribute these 
neigbhorhoods, more-or-less evenly over $\Gamma$, we select $J$ points 
$\{y_0,y_2,\ldots, y_{J-1}\}$ by first choosing $y_0 \in \Gamma$ arbitrarily, and subsequently incrementing $j$ and selecting $y_j$ to be the point of $\Gamma$ farthest from the previously selected points. 
Each selected point is the nexus of the $K$-nearest-neighborhood
\[
	\Omega_j = \NN_K(y_j).
\]
We choose $K=40$ and $J=10^4$, so that most of the $\Omega_j$ overlap with other neighborhoods, in the sense that they share points in $\Gamma$. In this case, each of the $\Omega_j$ is $0.03$-connected.

\begin{figure}[ht]
	\centering
	\includegraphics[scale=.4]{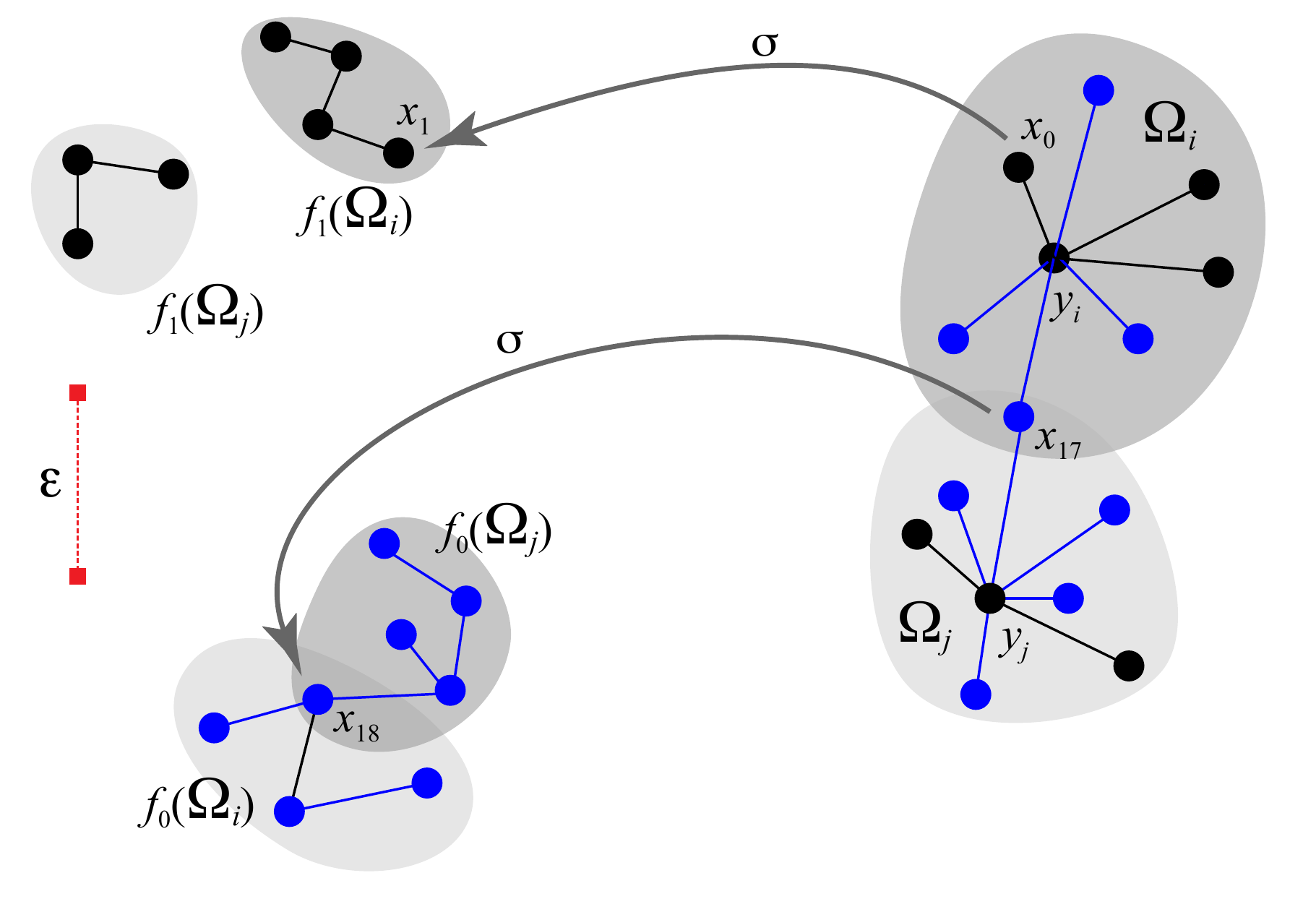}
	\caption{Separation of the time series $\Gamma$ into regimes. Here $\Omega_i$ and $\Omega_j$ represent the $8$-nearest neighborhoods of $y_i$ and $y_j$, respectively. They overlap, having $x_{17}$ in common. Each of the neighborhoods has two $\epsilon$-connected images under the shift $\sigma$. The pair that share $\sigma(x_{17}) = x_{18}$ are identified to be in the same regime, say $n=0$ so the $8$ points in this $\epsilon$-connected image set (and their preimages) are colored blue to indicate the common regime.
\label{fig:regimesketch}}
\end{figure}

The separation into regimes is accomplished as follows: whenever two ``overlapping" $\Omega_j$'s have $\epsilon$-connected image components that intersect, we identify them as being generated by the same $f_n$; see the sketch in \Fig{regimesketch}. More specifically, suppose that $\Omega_{j,k}$ is the set of points in $\Omega_j$ that generate the $k^{th}$ $\epsilon$-component of $\sigma(\Omega_j)$. These are distinguished by the following: whenever $\Omega_{j_1,k_1} \cap \Omega_{j_2,k_2} \neq \emptyset$, then the union of their images $\sigma(\Omega_{j_1,k_1}) \cup \sigma(\Omega_{j_2,k_2})$ will share a point as well, and thus be $\epsilon$-connected. In this case, the points in these images are selected as being generated by the same regime $f_n$. That is, $f_n(\Omega_{j_i,k_i}) = \sigma(\Omega_{j_i,k_i})$ for $i = 1,2$. \

The $\Omega_{j,k}$ can be thought of as nodes on an abstract graph. Whenever two of these neighborhoods share a point, an edge linking these nodes is added to the graph. Using this construction, the connected components of the resulting graph are selected as images of a fixed regime. Of course, we do not know which of the $f_n$'s is associated with which graph component unless we have prior knowledge of some values of the functions.

\begin{figure}[ht]
	\centering
	\subfigure{
		\includegraphics[scale=.40]{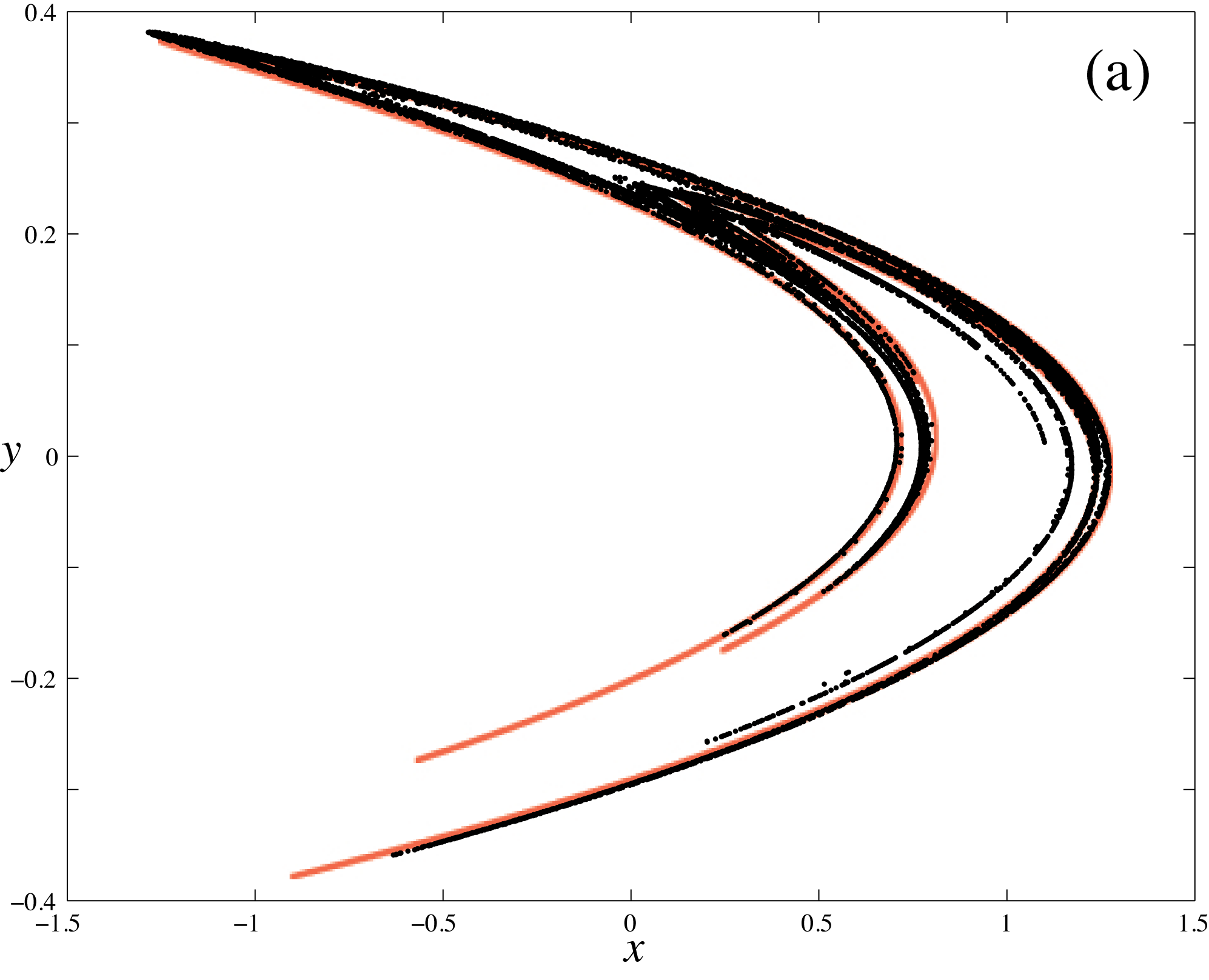}
	}
	\subfigure{
		\includegraphics[scale=.40]{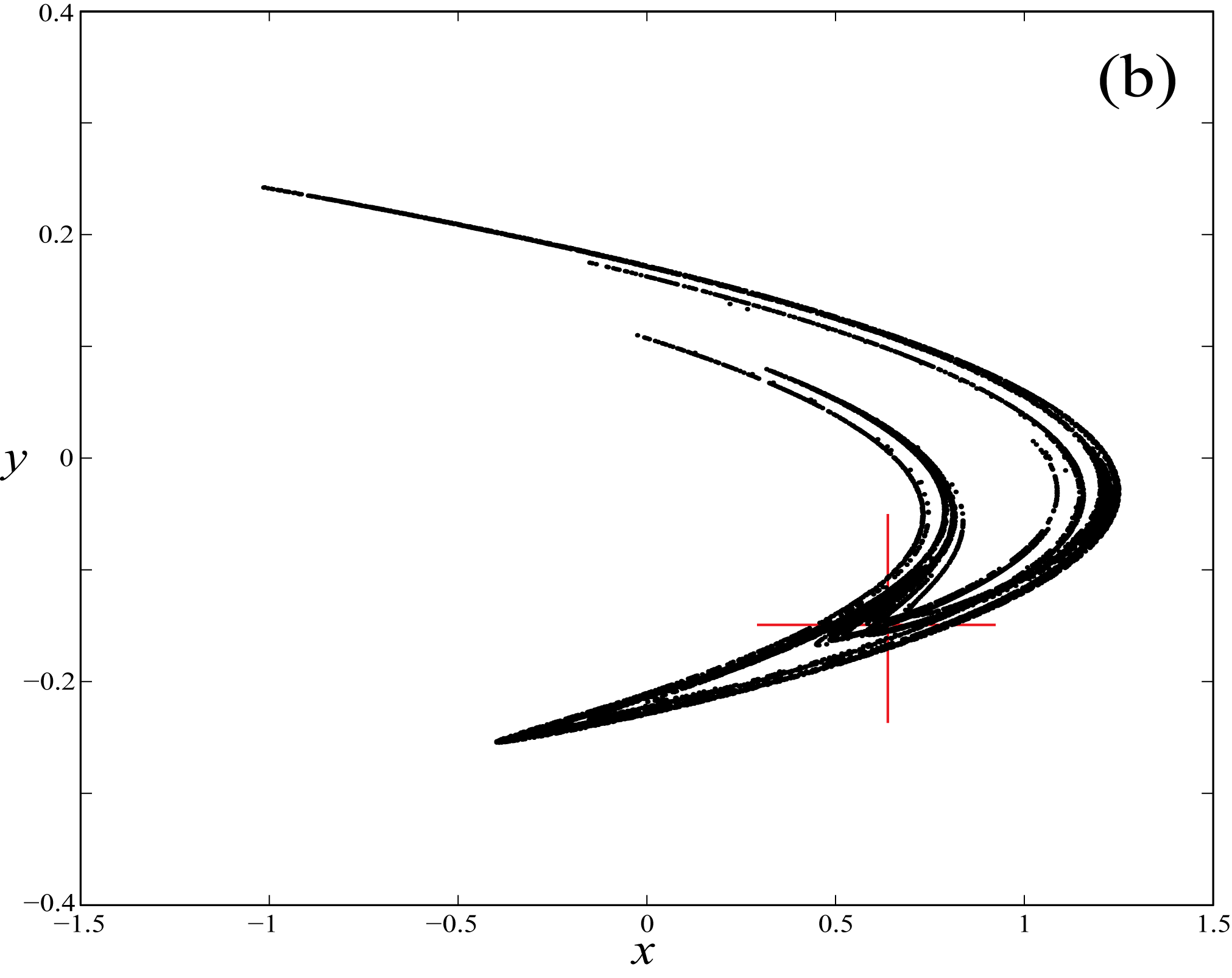}
	}
	\caption{Panel (a)  and panel (b) show the points identified as iterates of $f_0$ and $f_1$, respectively. These points can be approximately interpreted as a sampling of $f_0(\Gamma)$ ($f_1(\Gamma)$) Also shown, in grey (colored red, online) are points on the attractor of the H\'enon map $f_0$, and a cross at the position of the fixed point of $f_1$.\label{fig:separation}}
\end{figure}

For the trajectory of \Fig{ifs:trajectory} and using the covering by the $10^4$ neighborhoods $\Omega_j$, this algorithm generates two large connected graph components, one containing $14,724$ points and the other $14,815$ points. These sets of points are shown in the panels (a) and (b) of \Fig{separation}, respectively. Comparing these results with the known values of $n_t$ shows that every point in the first graph component has $n_t = 0$ and every point in the second has $n_t = 1$; that is, both the separation had \emph{no false positives}. There are an additional $465$ points of $\Gamma$ that are not in these two graph components. These unidentified points represent sparsely visited regions of the trajectory. 

It is no coincidence that the points identified to be images of $f_0$ in \Fig{separation}(a) appear to lie close to the attractor of the standard H\'enon map, which is shown in grey (light red, online) in the figure. Note, however, that even though the attractor for $f_1$ is a fixed point---the cross in the figure---the strong perturbation due to $f_0$ iterations causes the points in \Fig{separation}(b) to range far from the attractor of $f_1$.

\section{Computer Performance Dynamics} \label{sec:computer:performance:dynamics}
In this section, we describe the application of the regime separation algorithm to a time 
series obtained from an experimentally obtained computer performance analysis data set.  
A critical performance bottleneck in modern computer systems occurs in the efficient 
management of memory. The cache is the level of memory closest to the processor; it is 
preloaded with the data that the system \emph{thinks} it will need. When the system looks for 
a necessary piece of data in the cache and does not find it, it must load the data from 
main memory, resulting in a major performance slowdown. Such an event is called 
a \emph{cache miss}. 

\begin{figure}[th!]
	\centering
	\subfigure{
		\includegraphics[scale = .40]{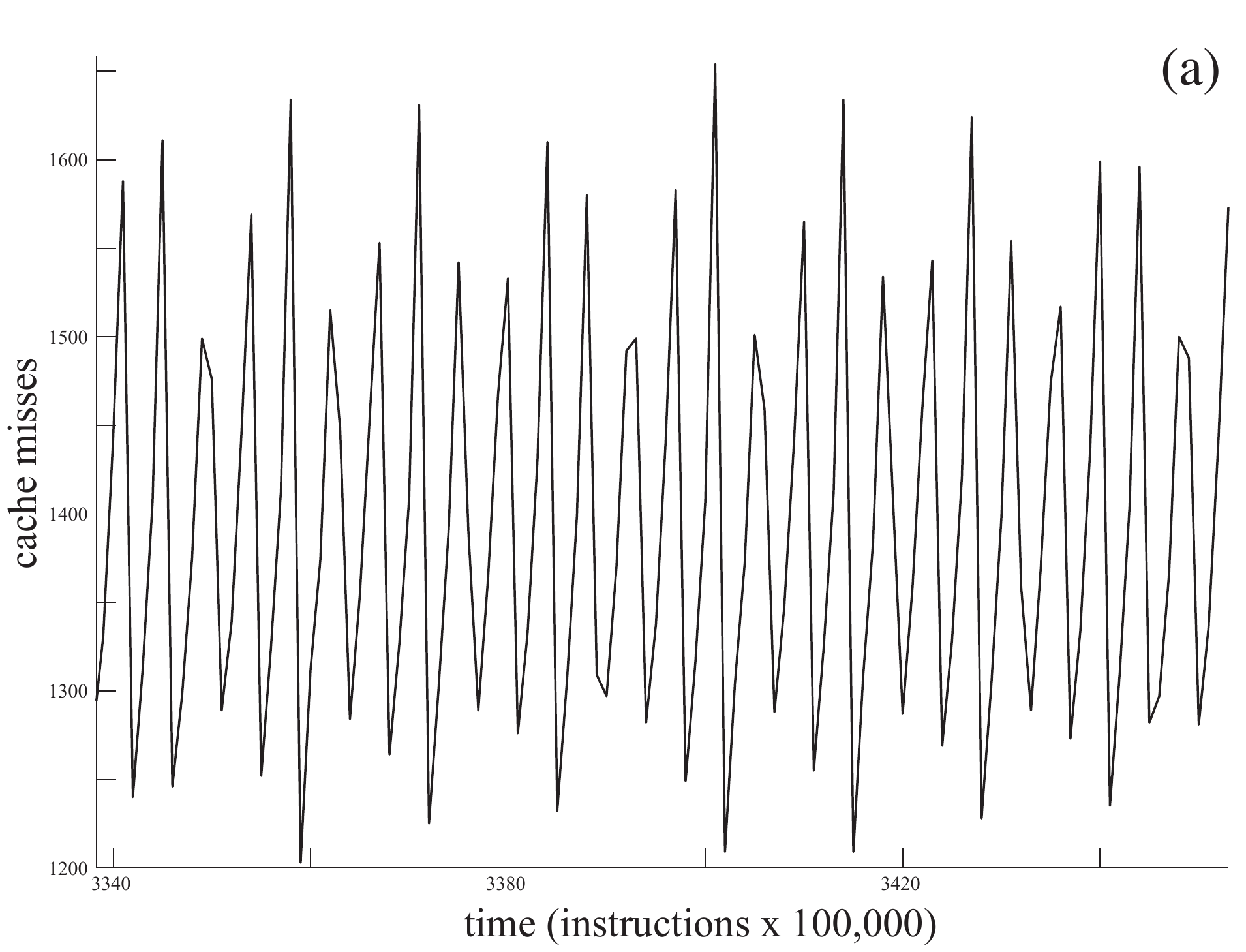}
		\label{fig:timedomain}
	}
	\subfigure{
		\includegraphics[scale = .4]{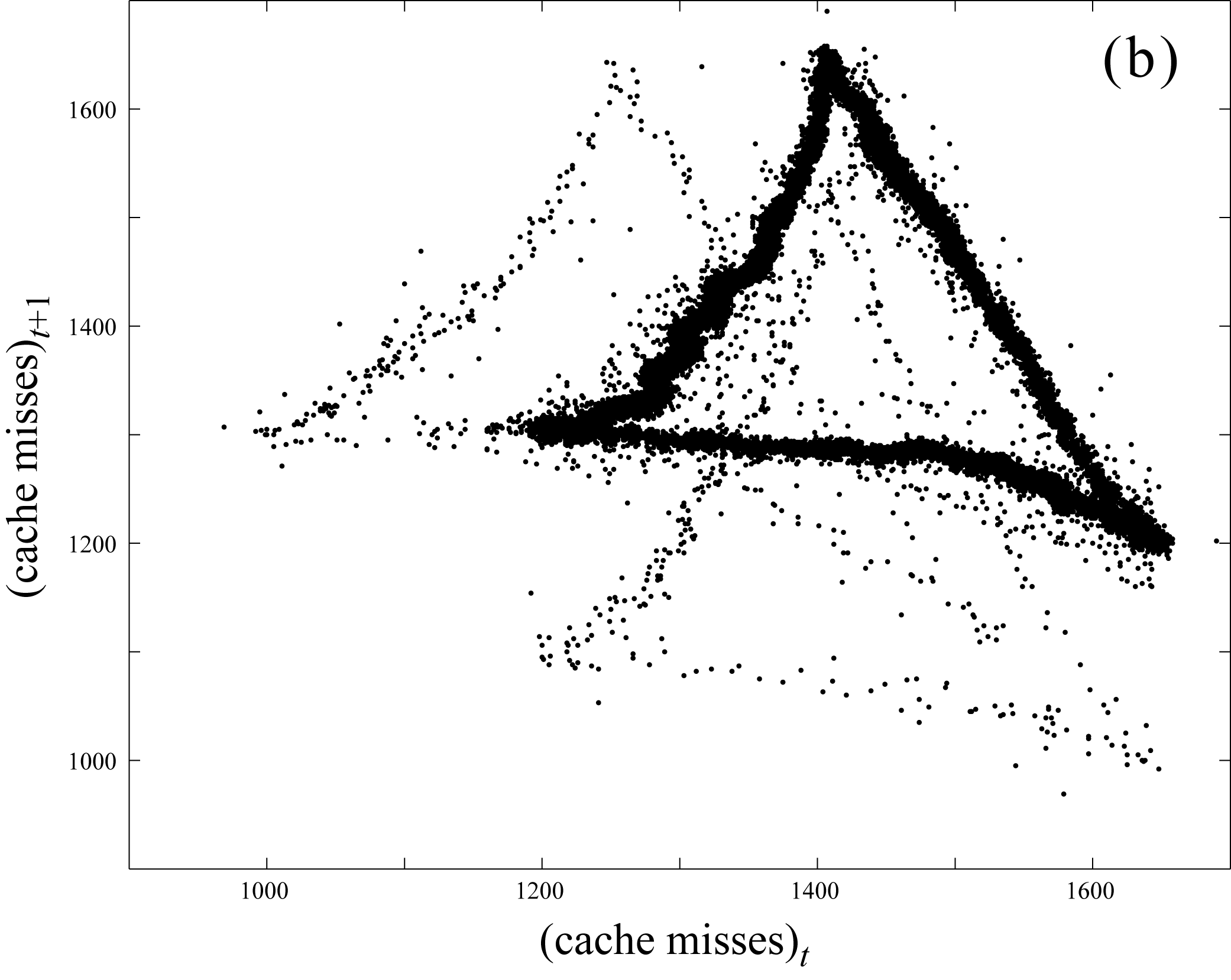}
		\label{fig:firstreturn}
	}
	\caption{(a) Cache misses per $10^5$ instructions observed during the execution of the code above.  (b) First return map of the time series from \Fig{timedomain} \label{fig:time:domain}}
\end{figure}

The experiment to investigate the frequency of cache misses consists of repeatedly running the simple C program:
\begin{lstlisting}[frame=tb] 
	for(i = 0; i < 255; i++) 
	  for(j = i; j < 255; j++) 
	    data[i][j] = 0; 
\end{lstlisting}
on an Intel Core2\textsuperscript{\textregistered} processor.
This code initializes the upper triangular portion of a matrix in
row-major order.  As the program runs, the hardware performance
monitors built into the processor monitor the memory
usage patterns---in particular, the rate of cache misses.
%
The program is interrupted every $10^5$ instructions and the
number of cache misses that occurred over that interval is recorded. We obtained a time series consisting of 86,107 points, from which we used a representative 60,000 point segment for the work presented in this paper. A snippet of the time series, along with the first return map of the 60,000 points selected is shown in \Fig{time:domain}. This data set has been studied previously and shown to exhibit chaotic
dynamics \cite{alexander2010measurement, mytkowicz2009computer}.

The time delay embedding process we use to analyze the computer performance trace involves the estimation of two parameters: the time
delay $\tau$  and the embedding dimension $m$. We first
choose $\tau$ using standard practices \cite{Sauer:1991lr,Takens:1981uq}. Based on the first minimum of the time-delayed average mutual information we choose $\tau = 10^5$ instructions. After choosing the delay $\tau$, the next step in the embedding
process is to estimate the dimension $m$. A standard strategy for this is the ``false nearest neighbor'' approach of Abarbanel and Kennel\cite{abarbanel1993local}.  
Based on this algorithm, we estimate that $10\leq m\leq 25$, then narrow
down that range to $m=12$ using other dynamical invariants \cite{mytkowicz2009computer}. Note that the first return map shown in \Fig{firstreturn}  is simply a two-dimensional projection of the 12-dimensional state space embedding, $\Gamma \subset \mathbb{R}^{12}$. For a more detailed discussion of these choices and our approach to estimating them see Garland and Bradley\cite{josh-IDA11} or Mytkowicz, et al.\cite{mytkowicz2009computer}.


The observation of the \emph{ghost triangles} in
\Fig{firstreturn}---seemingly reminiscent of three
overlapping attractors from an IFS---prompted us to apply our regime-separation algorithm to this data.  Because the two ghosts are much
more lightly sampled than the ``main'' triangle, our conjecture was
that the IFS consisted of three functions and that the switching
process prioritized one of the three. 

\begin{figure}
	\subfigure{
	\includegraphics[scale=.4,trim=.75in 0 0 0]{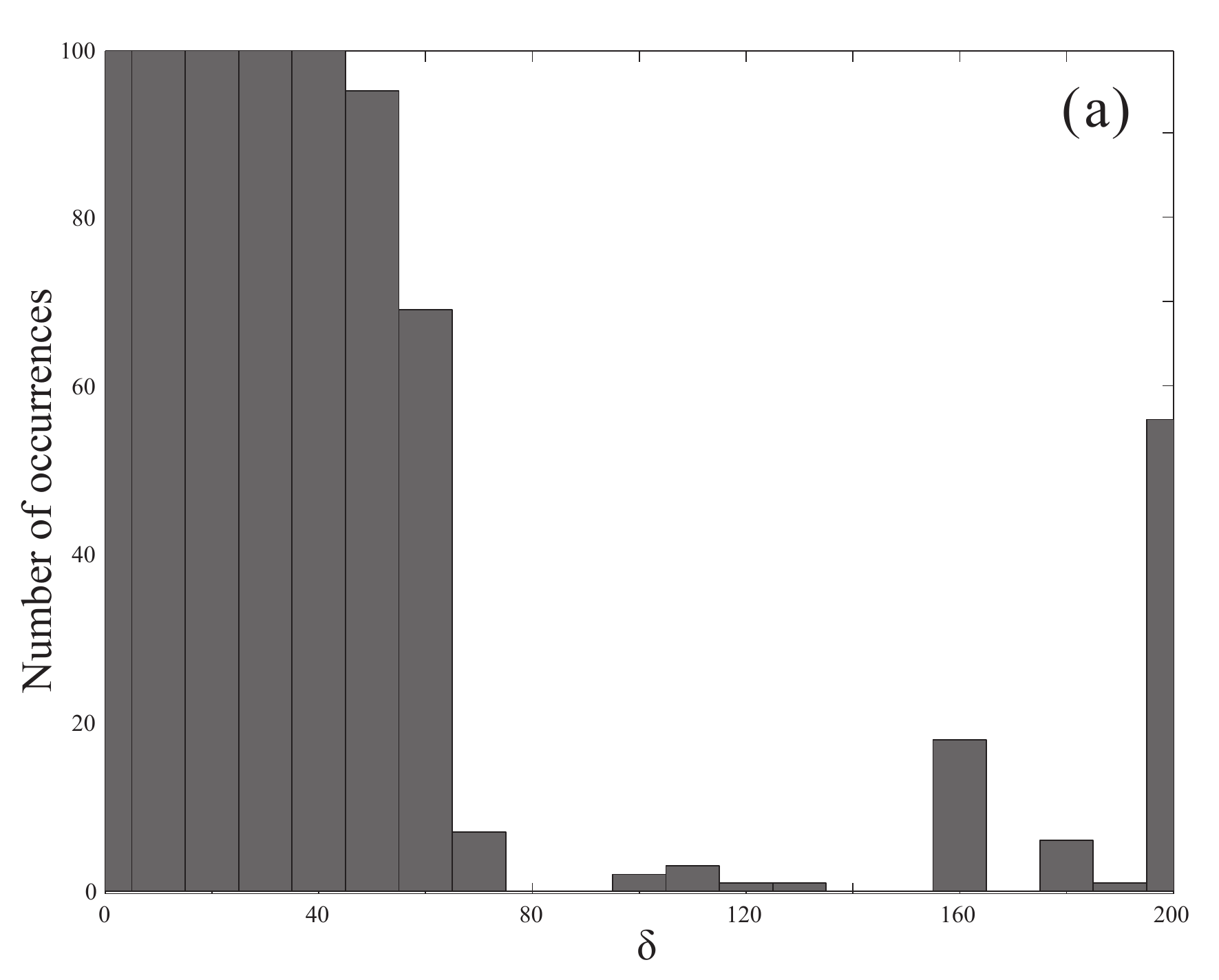}
	} \hspace{.3in}
	\subfigure{
	\includegraphics[scale=.4,trim=.75in 0 0 0]{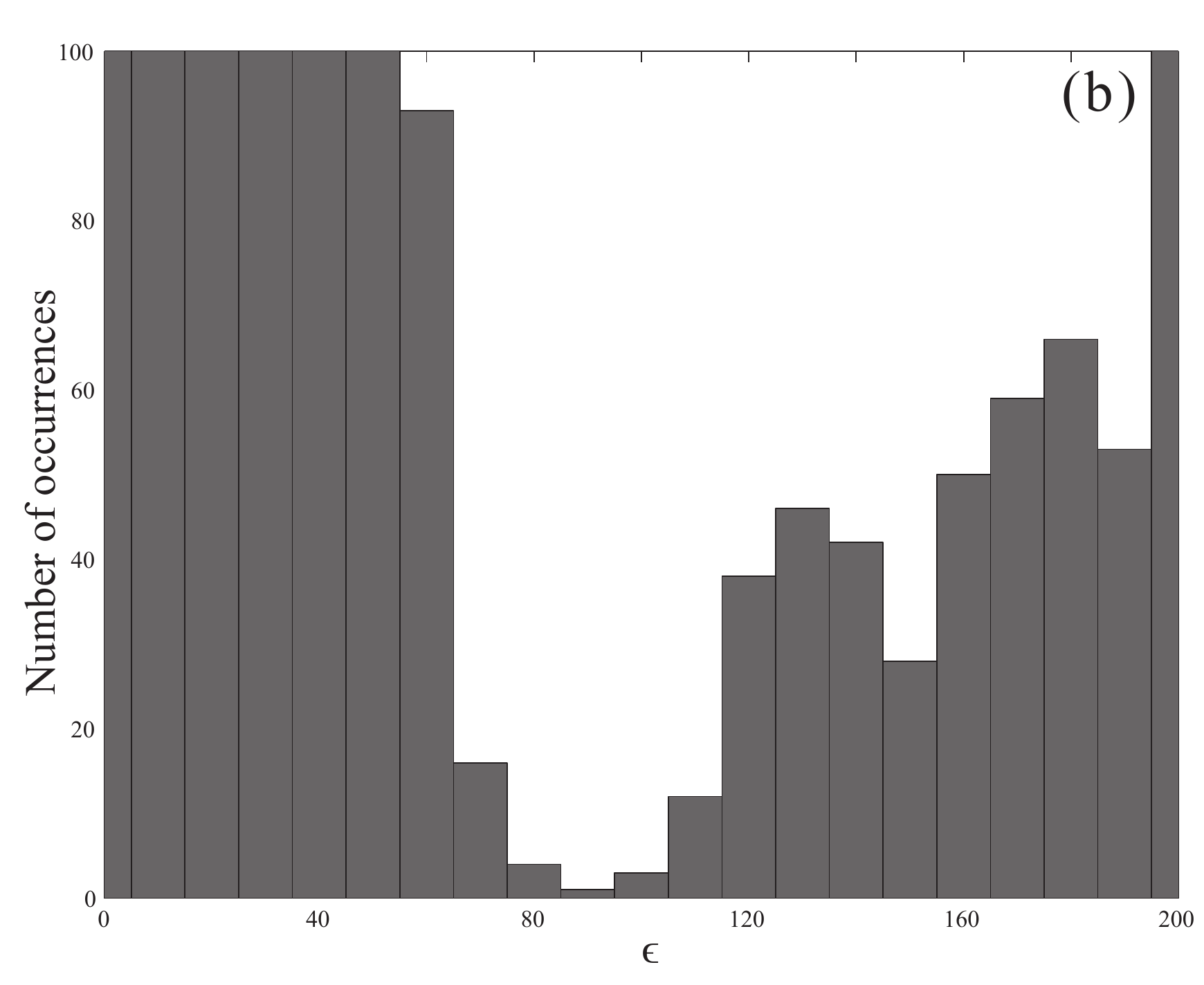}
	}
	\caption{Distance between each point on $\Gamma$---the 12-dimensional time-delay embedding of the computer performance data---and its nearest neighbor (a) and between the images of these two points (b).\label{fig:row}}
	
\end{figure}

For the analysis, we chose $\epsilon = 75$. This choice is justified by observation of the histograms shown in \Fig{row}. Figure \ref{fig:row}(a) is the histogram of the distances between each point in the time-series embedding and its nearest neighbor, i.e., the diameter of $\NN_2(x_j)$ in $\mathbb{R}^{12}$. Figure \ref{fig:row}(b), shows the histogram of distances between the same two points after iterating them both forward one time step, i.e., the diameter of $\sigma(\NN_2(x_j))$. The two `humps' indicate that, for a generic pair of nearest neighbors, the images of those neighbors are either within $\epsilon \approx 75$ of each other (iterates of the same $f_j$), or further than $\epsilon \approx 75$ apart (iterates of different $f_j$).


\begin{figure}
	\centering
	\includegraphics[scale=.4]{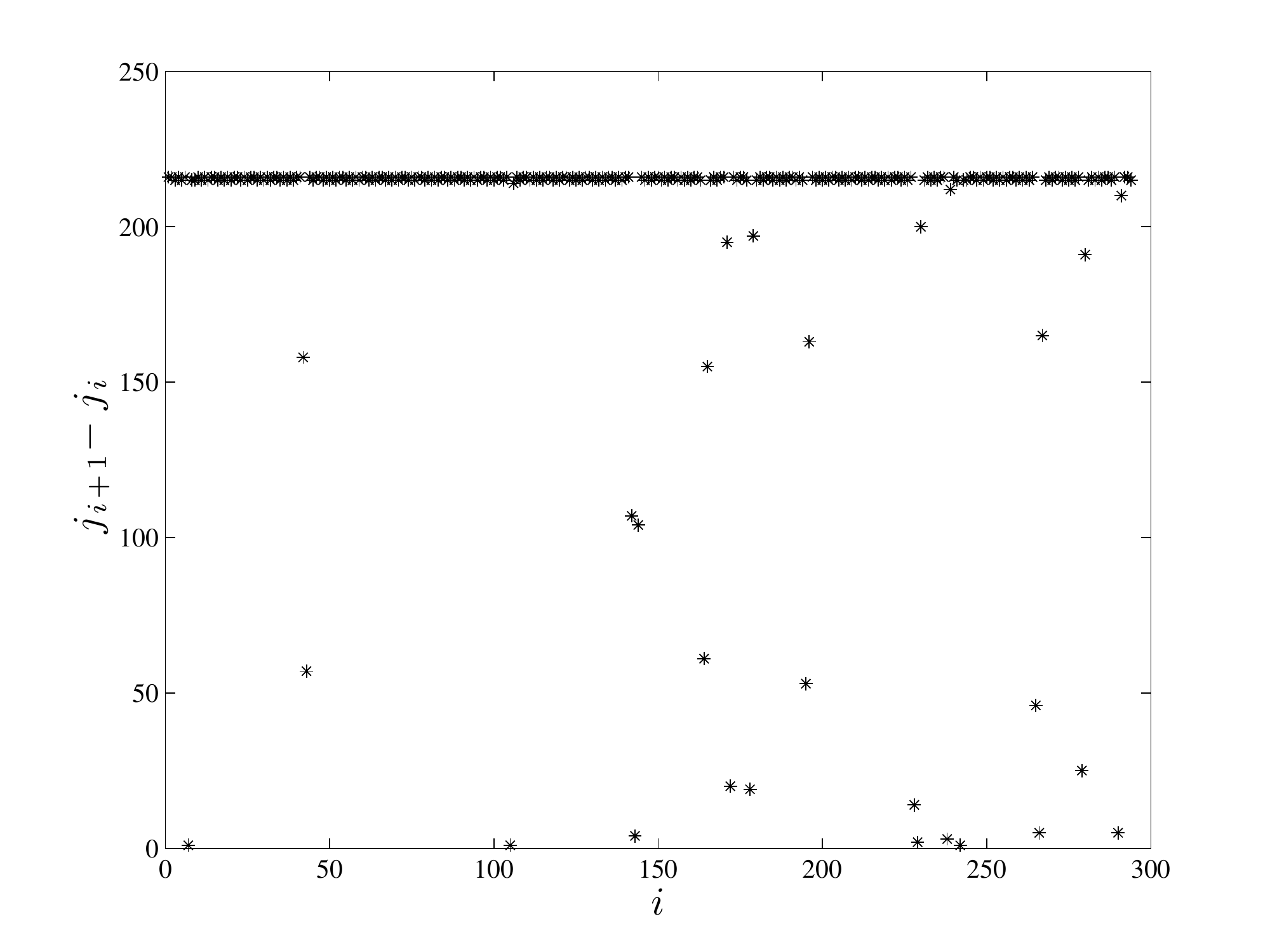}
	\caption{The first difference of the sequence $J = \{j_1,\ldots,j_K\}$ of ghost indices. That is $x_{j_i}$ is identified as a probably ghost point for $1 \le i \le K$. A point $(i,n)$ on the graph indicates that $j_{i+1} - j_{i} = n$. }
	\label{fig:ghostspacing}
\end{figure}

We separated the ghosts as follows: For each $x_j \in \Gamma$, we examine the image of $x_j$ and its nine nearest neighbors. That is, if $\sigma(\NN_{10}(x_j))$ consists of two $\epsilon$-connected components, then we identify the members of the smaller component (in cardinality) as candidates for points on the ghost triangle. Hence, we have a sequence $J = \{j_1,\ldots,j_K\}$ such that $x_{j_i} \in \Gamma$ is identified as a point on the ghost for $1 \le i \le K$. These points correspond to
the lower ghost of \Fig{firstreturn}; the second ghost is
just an image of the first---a necessary result of the symmetry
inherent in the time-delay embedding process. Furthermore, we note that the occurrence of ghost points is \emph{strongly} periodic, with a period of 215 measurements. This claim is motivated by the plot in \Fig{ghostspacing}. This plot is the first difference of the sequence $J$, i.e., the point $(i,n)$ indicates that the $(i+1)^{st}$ ghost occurs $n$ measurements after the $i^{th}$ ghost. Furthermore, the points on this graph that fall below the line $j_{i+1} - j_i = 215$ are the result of an `extra' identification in the middle of a period. For example, the first two such points have coordinates $(42,158)$ and $(43,57)$. Since $158 + 57 = 215$, one might hypothesize that the $42^{nd}$ ghost identification is actually spurious.

\begin{figure}[ht!]
	\centering
	\subfigure{
	\includegraphics[scale = .40]{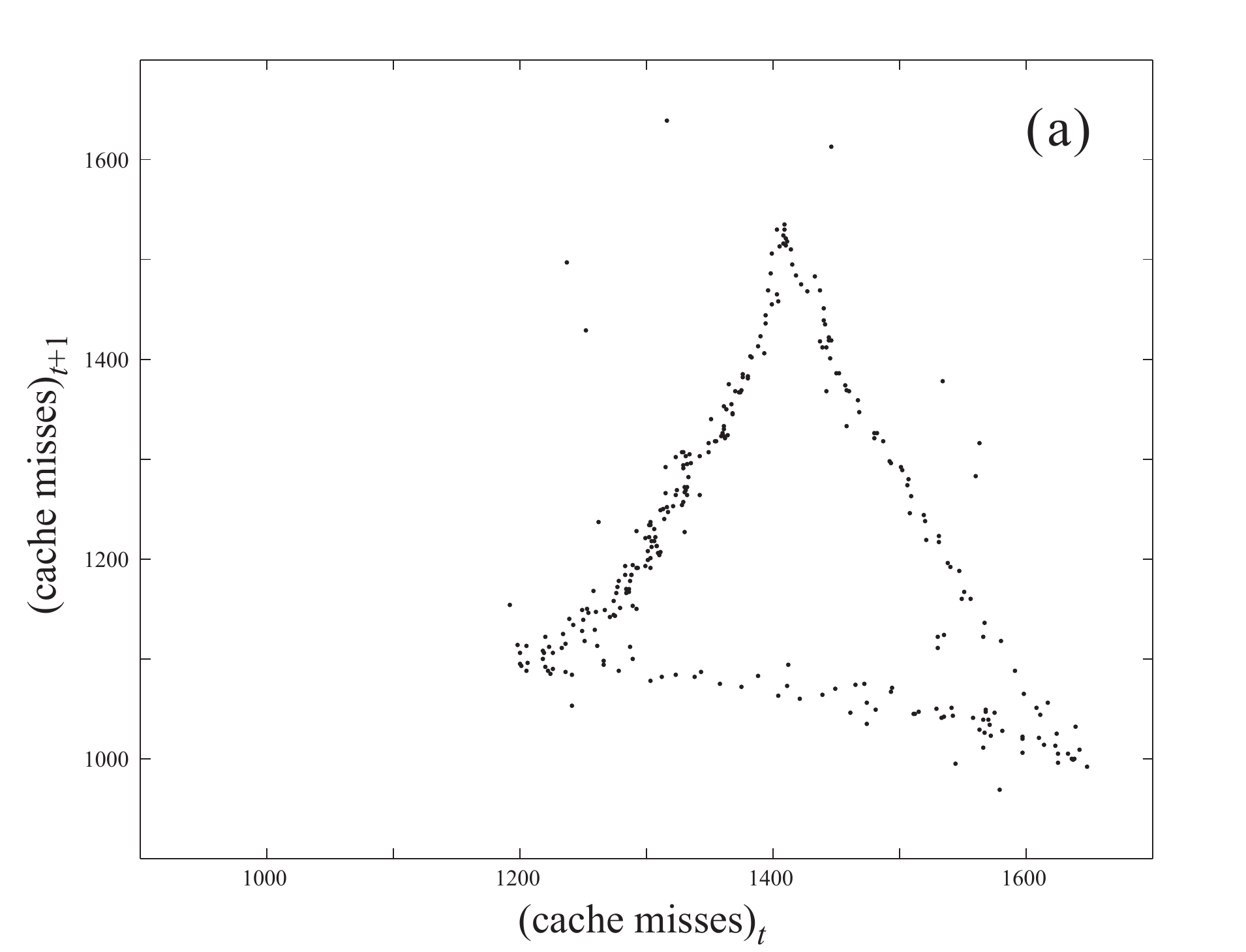}
	\label{fig:narrow}
	} \ \ \ \ 
	\subfigure{
	\includegraphics[scale = .40]{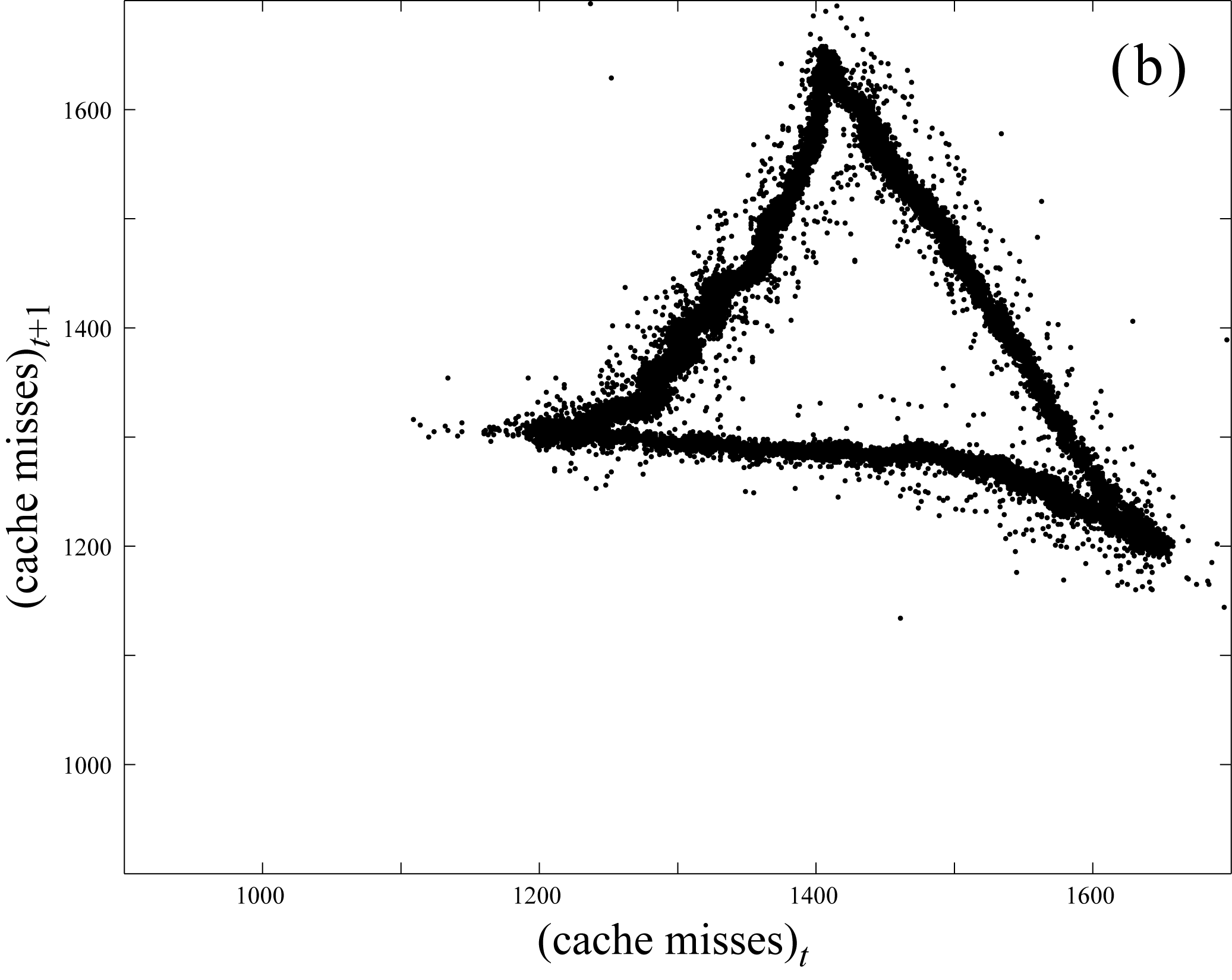}
	\label{fig:adjustment}
	}
	\caption{(a) The lower ghost triangle separated from the data of \Fig{firstreturn}. Each of the 295 points in this plot was identified as being $\epsilon$-disconnected ($\epsilon = 75$) from the images of an $\epsilon$-connected set of points. 
%
(b) A two dimensional time-delay embedding of the \emph{adjusted} time series obtained by adding 200 cache misses to the time-series values corresponding to each of the points from (a). \label{fig:ghost:replacement}}
\end{figure}

This analysis reinforces the hypothesis that there is a direct correspondence between the points on the ghosts and points in the time series that are periodically spaced by 215
measurements.  Moreover, each ghost point appears
to be shifted exactly 200 cache misses from the main triangle.  Indeed, adding $200$ cache misses to each of the points identified as parts of a ghost triangle, produces a time series that has the embedding shown in \Fig{adjustment}. Thus, for this case, not only is the regime identified, but the dynamics of the two components is shown to be simply related: by just a shift.

There is one issue remaining before we can model this data as an IFS: we only have access to \emph{measurements} of the states of the IFS.  That is, if
$X$ is the state space of the computer system and $\{f_0,\ldots,
f_k\}$ is a collection of maps on $X$, the observations correspond to the functions $\{h\circ f_0, \ldots,
h \circ f_k\}$ with a continuous measurement function $h:X \to \mathbb{R}$
that maps the state of the computer system to the number of cache
misses that occur over the given time interval. It can be shown that the functions $h\circ f_i$ are sufficient for studying topological and geometric properties of the $f_i$; a more-detailed treatment of the function $h$ can be found in Alexander et al.\cite{alexander2010measurement}.

Thus, if $h \circ f_0$ denotes the
dynamics associated with the main triangle, we can define $h \circ f_1(x)
= h \circ f_0(x) - 200$, and  $h \circ f_2(x) = h\circ f_0(y)$, where $y \in h^{-1}(h(x) + 200)$. 
The IFS consists of the state space $X$, the collection
$\{f_0,f_1,f_2\}$ of continuous maps, and the sequence
$\{n_j\} \subset \{0,1,2\}$, where
\[
	n_j = \left\{ \begin{array}{ll}
						1 \  & \ \ \text{if} \ j = 0\mod 215 \\
						2 \ & \ \ \text{if} \ j = 1\mod 215 \\
						0 \ & \ \ \text{otherwise} \\
				  \end{array} 
				\right.
\]
This model of the the cache-miss dynamics on the Intel Core2\textsuperscript{\textregistered} rests on the assumption that $f_1$ and $f_2$ can be described completely in terms of $f_0$. To verify this assumption, we
tested for determinism in the \emph{adjusted} dynamical system of
\Fig{narrow}.  We found that out of the 295 points so identified,
only 16 failed to lie in an $\epsilon$-connected image set in the
adjusted dynamical system. Consequently, $f_0$ appears to be a
continuous function and the IFS is an accurate model for this data
set.

Much of the usefulness of this model originates from the fact that we have
isolated the continuous function $f_0$.  In light of this, it is
reasonable to assume that $f_0$ is representative of some low-dimensional
dynamics that are present in the computer system, while $f_1$ and
$f_2$ represent a secondary piece of dynamics---in this case, perhaps
best described as `deterministic additive noise'.

\section{Conclusion} \label{sec:conclusion}
Many techniques for time-series analysis, such as in
Mischaikow et al.\cite{mischaikow1999construction}, explicitly require the time series
to be generated by a continuous function, and almost all of them
implicitly require that it be generated by a \emph{single} function.
For example, in Mytkowicz et al. \cite{mytkowicz2009computer}, time-series analysis of
the data of \Fig{timedomain} showed that it has a positive
Lyapunov exponent and fractal correlation dimension.  However, our
results show that that time series interleaves trajectories from
different dynamical systems---a property that can trip up traditional
time-series analysis techniques.  In the data studied in
Mytkowicz et al. \cite{mytkowicz2009computer}, this proved not to be an issue because a
single $f_0$ overwhelmingly dominated the dynamics.  When that is not
the case, problems can arise with traditional methods,
which are often formulated assuming the existence of long,
uninterrupted deterministic trajectories.  Some techniques, such as
those in Mischaikow et al.\cite{mischaikow1999construction}, do not explicitly require uninterrupted trajectories.  Using our topology-based approach, one could pull apart the time-series data into individual regimes and study the dynamics of each of the
$f_i$ independently.



In conclusion, we have described an algorithm for detection and separation of a signal that is generated by continuous, deterministic dynamics punctuated by regime shifts. The algorithm handles shifts that result from stochastic or deterministic processes: it applies whenever the dynamics are described by an iterated function system. Time-series data from a computer performance analysis experiment were shown to fit this model. More generally, we claim that iterated
function systems are a natural model for complex computer programs, which---we hypothesize---have regime shifts as their execution moves through different parts of the code. 

IFS models provide a natural framework for data analysis in a wide range of fields: whenever the physical system generating the data is characterized by continuous systems, that  are punctuated by discontinuous regime shifts. 
Another area in which the regime separation technique is particularly
appropriate is digital communication channels. Indeed, (hyperbolic)
iterated function systems are known to provide useful models of these
channels \cite{broomhead2004iterated}. A channel corresponds to an
electrical circuit externally driven by a digital signal, and the
discrete input signal corresponds to the regime sequence. Thus, the
behavior of the circuit corresponds to the actions of a discrete set of continuous
dynamical systems. A fundamental problem in this context is
\emph{channel equalization}, the reversal of distortion that is incurred by
transmission through a channel. This is precisely the determination of
the input signal sequence from a sequence of output values---i.e.,
regime separation. Channel equalization is straightforward for linear
dynamics because the IFS attractors in these situations tend to be
non-overlapping. However, more-realistic, nonlinear IFS models have
overlapping attractors. We believe that our methods can be
successfully used for channel equalization in this context.

Challenges that remain to be addressed include finding an efficient
implementation for high-dimensional data and dealing with systems that
have traditional (e.g. Gaussian) noise in addition to regime shifts. Furthermore, we have not addressed the nature of the switching process itself. Once the regime shifts have been determined, the next natural question to ask is whether or not the switching is deterministic or stochastic, and if one can determine the rule for switching between regimes.

\end{document}